\documentclass[12pt]{amsart}
\usepackage{a4wide,enumerate,color,graphicx}
\usepackage{amsmath}

 \newtheorem{definition}{Definition}
 \newtheorem{thm}{Theorem}
 
 \newtheorem{lem}[thm]{Lemma}
 
 \theoremstyle{definition}
 
 \theoremstyle{remark}

 \DeclareMathOperator{\dive}{div}

\begin{document}
\title[]{Bounded weak solutions to cross-diffusion semiconductor model with electron-hole scattering}

\author[]{Xi Lin}
\address{{Department of Mathematics and Physics, Guangzhou Maritime University,
Guangzhou 510725, Guangdong Province, China}}
\email{linxi@gzmtu.edu.cn}

\address{}
\email{}
%\thanks{The author acknowledges support from the National Natural Science
%Foundation of China (NSFC), grant 12471206.}

%\date{\today}

%\tableofcontents

%\newpage

%\begin{abstract}

%\end{abstract}
%\maketitle

\begin{abstract}
Semiconductor model is a system of parabolic partial differential equations with cross-diffusion phenomenon. Previous results showed that a weak solution exists and is not bounded in general. So semiconductor model was categorized as a cross-diffusion system without bounded weak solutions. In this work, we show that once the initial value is bounded, there exists a weak solution that is also bounded. The entropy method is a major tool in global existence analysis of cross-diffusion systems. We notice that traditional entropies in volume-filling cases may not provide required positive semi-definiteness result for the existence proof. In this situation, a transformation of variables technique has been applied. The product between Hessian matrix of the entropy and replacement diffusion matrix is positive semi-definite, then we apply the entropy method to show semiconductor model has a bounded weak solution.
\end{abstract}
\maketitle

\section{Introduction}

Cross-diffusion systems are parabolic partial differential equations. They arise naturally in multicomponent systems. This work is a continuation of articles \cite{1,2,3}. In these references, the authors investigated analytical properties of global-in-time solutions $u=(u_1,...,u_n):\Omega\times(0,\infty)\to\mathbb{R}^n$ to a parabolic system 
\begin{equation}\label{equation}
\partial_tu_i=\dive\Big(\sum_{j=1}^{n}A_{ij}(u)\nabla u_j\Big)\text{ in }Q_T=\Omega\times(0,T),\text{ }T>0,\text{ }1\leq i\leq n,
\end{equation}
where $\Omega\subset\mathbb{R}^{d}$, $d\geq 1$ is a bounded domain. Equations subject to initial and no-flux boundary conditions
\begin{equation}\label{initial}
\sum_{j=1}^{n}A_{ij}(u)\nabla u_j\cdot \nu=0 \text{ }\text{   on  } \partial{\Omega},\text{ } t>0,\text{ } u_i(0)=u^{0}_i \text{  in  } \Omega,\text{ }1\leq i\leq n,
\end{equation}
where $\nu$ is an exterior unit normal vector to $\partial\Omega$ and $u^{0}_i$ is the initial datum. The system is (\ref{equation})-(\ref{initial}) with $n=2$, $u=(u_1,u_2)$, and the diffusion matrix $A(u)$ is given by
\begin{equation}\label{diffusiontwo}
\begin{aligned}
A(u)&=\frac{1}{1+\mu_2u_1+\mu_1u_2}\left(\begin{array}{cccccc}
\mu_1(1+\mu_2u_1) & \mu_1\mu_2u_1 \\
 \mu_1\mu_2u_2& \mu_2(1+\mu_1u_2) \\
\end{array}\right).
\end{aligned}
\end{equation}

The system of equations (\ref{equation})-(\ref{diffusiontwo}) is called semiconductor model with electron-hole scattering, which is listed in page 1976, \cite{2}. The model has its origination in \cite{4,5}.

Semiconductor model describes how electron-hole scattering affects the carrier transportation in high-injection situations. In this model, $u_1$ refers to the electron-hole density, $u_2$ refers to the hole density, $\mu_1$ and $\mu_2$ are positive mobility constants for electrons and holes, respectively. 

The existence of a global weak solution to (\ref{equation})-(\ref{diffusiontwo}) was shown in \cite{1}, and this solution is not bounded in general. It was shown in \cite{3} that when $\mu_1=\mu_2=1$, this weak solution is unique.

Cross-diffusion systems are divided into two groups in \cite{2}. Let us denote $\mathcal{D}\subset\mathbb{R}^n$ a bounded open set and $\overline{\mathcal{D}}$ is the closure set of $\mathcal{D}$. Once the weak solution $u$ of (\ref{equation})-(\ref{initial}) satisfies that $u\in\overline{\mathcal{D}}$, then we call $u$ is bounded. Otherwise $u$ is not a bounded weak solution. 

Degenerate systems \cite{6,7} and Maxwell-Stefan equations \cite{8,9} are typical cross-diffusion systems with bounded weak solutions. The Shigesada-Kawasaki-Teramoto type population system \cite{10,11} is a cross-diffusion system with its weak solutions not bounded in general. Semiconductor model was categorized as a model with unbounded weak solutions in \cite{2}.

In this work, we show that once the initial value is bounded, then (\ref{equation})-(\ref{diffusiontwo}) has a bounded weak solution.
Once the initial value is not bounded, then a weak solution exists and is not bounded. To the best of our knowledge, few cross-diffusion systems possess such properties. Authors of previous references were not aware of the fact that semiconductor model actually is also a cross-diffusion system with bounded weak solutions.

The original equations for semiconductor model in \cite{1} are
\begin{equation}\nonumber
\begin{aligned}
\partial_tp-\dive J_p=R(p,q),\quad\partial_tq-\dive J_q=R(p,q),\quad R(p,q)=\frac{p^2_i-pq}{\tau_0+\tau_pp+\tau_qq},
\end{aligned}
\end{equation}
where $p$ and $q$ are electron and hole densities of the semiconductor crystal, $p_i,\tau_0,\tau_p,\tau_q$ are positive constants, and
\begin{equation}\nonumber
\begin{aligned}
&J_p=\mu_{pp}(p,q)(\nabla p-p\nabla\psi)+\mu_{pq}(p,q)(\nabla q+q\nabla\psi),\\&J_q=\mu_{qp}(p,q)(\nabla p-p\nabla\psi)+\mu_{qq}(p,q)(\nabla q+q\nabla\psi),\quad\lambda^2\Delta\psi=p-q-C(x).
\end{aligned}
\end{equation}

Mobilities $\mu_{pp},\mu_{pq},\mu_{qp},\mu_{qq}$ are defined as
\begin{equation}\nonumber
\begin{aligned}
&\mu_{pp}(p,q)=\frac{\mu_p(1+\mu_qp)}{1+\mu_qp+\mu_pq},\quad\mu_{pq}(p,q)=\frac{\mu_p\mu_qp}{1+\mu_qp+\mu_pq},\\&\mu_{qq}(p,q)=\frac{\mu_q(1+\mu_pq)}{1+\mu_qp+\mu_pq},\quad\mu_{qp}(p,q)=\frac{\mu_q\mu_pq}{1+\mu_qp+\mu_pq},
\end{aligned}
\end{equation}
and $\mu_p$, $\mu_q$ are positive constants. The function $\psi$ is the electrostatic potential, $C(x)$ is the so-called doping profile. $R(p,q)$ is called the Shockley-Read-Hall term, which models generation-recombination process. In (\ref{equation})-(\ref{diffusiontwo}), the electrostatic potential $\psi=0$, and the Shockley-Read-Hall term $R(p,q)=0$. Notations $p,q,\mu_p,\mu_q$ are replaced correspondingly by $u_1,u_2,\mu_1,\mu_2$.

In \cite{1}, the regularity property of the weak solution $u$ is rather weak. In this manuscript, the definition of weak solutions to (\ref{equation})-(\ref{diffusiontwo}) has been modified and given as Definition \ref{defweak}. In the following discussion, we denote $\mathcal{D}_M=\{u=(u_1,u_2):u_1+u_2<M,u_1,u_2>0,M>0\}$, and $\overline{\mathcal{D}}_{M^{\prime}}$ is the closure of $\mathcal{D}_{M^{\prime}}=\{u=(u_1,u_2):u_1+u_2<M^{\prime},u_1,u_2>0,M^{\prime}>0\}$.
\begin{definition}\label{defweak}
For every $T>0$, $u^0=(u^0_1,u^0_2)\in\mathcal{D}_M$, if for all test functions $\phi=(\phi_1,\phi_2)$, $\phi_1,\phi_2\in L^2(0,T;H^1(\Omega))$, $u=(u_1,u_2)\in\overline{\mathcal{D}}_{M^{\prime}}$, and satisfies that
\begin{equation}\nonumber
\begin{aligned}
&\int_{0}^{T}\langle\partial_tu_1,\phi_1\rangle dt+\int_{Q_T}\frac{\mu_1(1+\mu_2u_1)\nabla u_1+\mu_1\mu_2u_1\nabla u_2}{1+\mu_2u_1+\mu_1u_2}\cdot\nabla\phi_1dxdt=0,\\&\int_{0}^{T}\langle\partial_tu_2,\phi_2\rangle dt+\int_{Q_T}\frac{\mu_2(1+\mu_1u_2)\nabla u_2+\mu_1\mu_2u_2\nabla u_1}{1+\mu_2u_1+\mu_1u_2}\cdot\nabla\phi_2dxdt=0,
\end{aligned}
\end{equation}
with the property that $u_i\in L^2(0,T;H^1(\Omega))$, $\partial_tu_i\in L^2(0,T;H^1(\Omega)^{\prime})$, $i=1,2$, then we call $u$ is a bounded weak solution of $(\ref{equation})$-$(\ref{diffusiontwo})$. In here, $M^{\prime}>M$ is a constant.
\end{definition}

We try to show Theorem \ref{mainresult} in following sections, which is the main result of this work.
\begin{thm}\label{mainresult}
Let $u^0=(u^0_1,u^0_2)$ be such that $u^0\in\mathcal{D}_M$. There exists a weak solution $u\in\overline{\mathcal{D}}_{M^{\prime}}$ of $(\ref{equation})$-$(\ref{diffusiontwo})$, with $M^{\prime}>M$ a constant.
\end{thm}

In \cite{3}, a key step of the proof was to show that semiconductor model has a bounded weak solution, once $\mu_1=\mu_2$ and the initial value is bounded. For $\mu_1\neq\mu_2$ case, when the initial value is bounded, whether a bounded weak solution exists or not remains unknown.

It was conjectured in \cite{3} that if $\mu_1\neq\mu_2$, the weak solution of (\ref{equation})-(\ref{diffusiontwo}) is still unique. This is one of central topics in the study of parabolic partial differential equations. After this work, we are left to consider uniqueness of these bounded weak solutions. It is sufficient to show that once the initial value is bounded, the bouneded weak solution is unique. Then by the method in \cite{3}, we can show that the weak solution of semiconductor model is unique.

In \cite{12,13,14,15,16}, the authors discussed parabolic systems of the form (\ref{equation})-(\ref{initial}). The existence of bounded weak solutions to semiconductor model is a continuation of these discussions. Semiconductor model is a typical cross-diffusion system with strong physical background, so it was investigated at a rather early time. In this work, we realize that semiconductor model is also a rather special cross-diffusion system. Classification of cross-diffusion systems criterion in \cite{2} may not be appropriate for semiconductor model, as this model belongs to both categories.

\section{Major Difficulties}

In global existence analysis of cross-diffusion systems, we mainly rely on the entropy method to show weak solutions exist. A major step is to find an appropriate entropy $h(u)$, such that $h^{\prime\prime}(u)A(u)$ is positive semi-definite. In here, $H(u)=h^{\prime\prime}(u)$ is the Hessian matrix of the entropy density $h(u)$. For every $1\leq i,j\leq n$, $H_{ij}(u)=\partial^2h/\partial u_i\partial u_j$. Positive semi-definite refers to the generalized definition: $h^{\prime\prime}(u)A(u)$ satisfies that for every $z=(z_1,...,z_n)^{\mathrm{T}}\in\mathbb{R}^n$, $z^{\mathrm{T}}h^{\prime\prime}(u)A(u)z\geq 0$, but $h^{\prime\prime}(u)A(u)$ is not symmetric in general.

If such an entropy $h(u)$ exists, then a systematic proof scheme can be applied to show a weak solution exists to a system of cross-diffusion parabolic partial differential equations. For example, in page 1976, \cite{2}, the entropy 
\begin{equation}\nonumber
\begin{aligned}
h_1(u)=u_1(\log u_1-1)+u_2(\log u_2-1)
\end{aligned}
\end{equation}
was applied, which follows
\begin{equation}\nonumber
\begin{aligned}
h^{\prime\prime}_1(u)A(u)&=\frac{1}{1+\mu_2u_1+\mu_1u_2}\left(\begin{array}{cccccc}
\frac{1}{u_1} & 0 \\
0 & \frac{1}{u_2} \\
\end{array}\right)\cdot\left(\begin{array}{cccccc}
\mu_1(1+\mu_2u_1) & \mu_1\mu_2u_1 \\
 \mu_1\mu_2u_2& \mu_2(1+\mu_1u_2) \\
\end{array}\right)\\&=\frac{1}{1+\mu_2u_1+\mu_1u_2}\left(\begin{array}{cccccc}
\frac{\mu_1}{u_1}+\mu_1\mu_2 & \mu_1\mu_2 \\
\mu_1\mu_2 & \frac{\mu_2}{u_2}+\mu_1\mu_2 \\
\end{array}\right).
\end{aligned}
\end{equation}

For all $z=(z_1,z_2)^{\mathrm{T}}\in\mathbb{R}^2$, $u_1,u_2>0$, we have
\begin{equation}\nonumber
\begin{aligned}
z^{\mathrm{T}}h^{\prime\prime}_1(u)A(u)z=\frac{1}{1+\mu_2u_1+\mu_1u_2}\big(\frac{\mu_1}{u_1}z^2_1+\frac{\mu_2}{u_2}z^2_2+\mu_1\mu_2(z_1+z_2)^2\big)\geq 0.
\end{aligned}
\end{equation}

By applying the entropy $h_1(u)$, we can show that a weak solution of (\ref{equation})-(\ref{diffusiontwo}) exists. Since the domain for $h_1(u)$ is $\{u=(u_1,u_2):u_1,u_2>0\}$, whether this weak solution is bounded or not remains unknown, even if the initial value is bounded. 

We mentioned that Maxwell-Stefan equations and degenerate system are important cross-diffusion systems with bouneded weak solutions. In global existence proof of Maxwell-Stefan equations \cite{8}, the entropy was
\begin{equation}\label{maxwellentropy}
\begin{aligned}
h(\rho^{\prime})=k\sum_{i=1}^{N+1}x_i(\ln x_i-1)+k,\quad\rho^{\prime}=(\rho_1,...,\rho_N).
\end{aligned}
\end{equation}
 
In (\ref{maxwellentropy}), $x_i$ and $\rho_i$ are related by $x_i=\rho_i/(kM_i)$, $\rho_i>0$, with $k=\sum_{i=1}^{N+1}(\rho_i/M_i)$, and $\sum_{i=1}^{N+1}\rho_i=\sum_{i=1}^{N+1}x_i=1$. $M_i>0$, $i=1,...,N+1$ are constants. The domain for $h(\rho^{\prime})$ in (\ref{maxwellentropy}) is $\{\rho^{\prime}=(\rho_1,...,\rho_N):\sum_{i=1}^{N}\rho_i<1$, $\rho_i>0$\}.

For degenerate system \cite{7}, the entropy was
\begin{equation}\label{degenerateentropy}
\begin{aligned}
h(u)=\sum_{i=1}^{n}(u_i\log u_i-u_i+1)+\int_{a}^{u_{n+1}}\log q(s)ds+\chi(u).
\end{aligned}
\end{equation}

In (\ref{degenerateentropy}), $a\in(0,1]$ is a constant, $\sum_{i=1}^{n+1}u_i=1$, $u_i>0$, and $q(s),\chi(u)$ are functions defined in (7)-(8), \cite{7}. The domain for $h(u)$ in (\ref{degenerateentropy}) is $\{u=(u_1,...,u_n):\sum_{i=1}^{n}u_i<1,u_i>0\}$.

For semiconductor model, we want to show a weak solution $u\in\overline{\mathcal{D}}_M$ exists. Let us denote $u_3=M-u_1-u_2$, and choose the entropy
\begin{equation}\label{h2u}
\begin{aligned}
h_2(u)=u_1(\log u_1-1)+u_2(\log u_2-1)+u_3(\log u_3-1).
\end{aligned}
\end{equation}

We notice that entropies in (\ref{maxwellentropy})-(\ref{degenerateentropy}) are intimately related to $h_2(u)$ in (\ref{h2u}). $h_2(u)$ is a major functional in showing a bounded weak solution exists to a cross-diffusion system. For $u\in\mathcal{D}_M$, we have
\begin{equation}\nonumber
\begin{aligned}
&(1+\mu_2u_1+\mu_1u_2)h^{\prime\prime}_2(u)A(u)\\&=\Big\{\left(\begin{array}{cccccc}
\frac{1}{u_1} & 0 \\
 0 & \frac{1}{u_2} \\
\end{array}\right)+\frac{1}{u_3}\left(\begin{array}{cccccc}
1 & 1 \\
1 & 1 \\
\end{array}\right)\Big\}\cdot\left(\begin{array}{cccccc}
\mu_1(1+\mu_2u_1) & \mu_1\mu_2u_1 \\
 \mu_1\mu_2u_2 & \mu_2(1+\mu_1u_2) \\
\end{array}\right)\\&=\left(\begin{array}{cccccc}
\frac{\mu_1}{u_1} & 0 \\
 0 & \frac{\mu_2}{u_2} \\
\end{array}\right)+\frac{M\mu_1\mu_2}{u_3}\left(\begin{array}{cccccc}
1 & 1 \\
1 & 1 \\
\end{array}\right)+\frac{1}{u_3}\left(\begin{array}{cccccc}
\mu_1 & \mu_2 \\
\mu_1 & \mu_2 \\
\end{array}\right),
\end{aligned}
\end{equation}
which follows for $z=(z_1,z_2)^{\mathrm{T}}\in\mathbb{R}^2$,
\begin{equation}\nonumber
\begin{aligned}
&(1+\mu_2u_1+\mu_1u_2)\cdot z^{\mathrm{T}}h^{\prime\prime}_2(u)A(u)z\\&=\frac{\mu_1}{u_1}z^2_1+\frac{\mu_2}{u_2}z^2_2+\frac{M\mu_1\mu_2}{u_3}(z_1+z_2)^2+\frac{1}{u_3}(\mu_1z^2_1+(\mu_1+\mu_2)z_1z_2+\mu_2z^2_2).
\end{aligned}
\end{equation}

If $\mu_1\neq\mu_2$, it does not hold that for every $z_1,z_2\in\mathbb{R}$, $$\mu_1z^2_1+(\mu_1+\mu_2)z_1z_2+\mu_2z^2_2\geq 0,$$ then we cannot show that for $z=(z_1,z_2)^{\mathrm{T}}\in\mathbb{R}^2$, $u\in\mathcal{D}_M$, $z^{\mathrm{T}}h^{\prime\prime}_2(u)A(u)z\geq 0$.

In above computations, we notice that $h^{\prime\prime}_2(u)A(u)$ is not positive semi-definite, though $h^{\prime\prime}_1(u)A(u)$ is positive semi-definite. A weak solution that is not bounded can be shown exists, but no existence result can be derived for a bounded weak solution, if the initial value is bounded. So the existence of a bounded weak solution for semiconductor model was not proved in previous references.

Actually, the product between Hessian matrix of the entropy and diffusion matrix is not positive semi-definite in general. For Shigesada-Kawasaki-Teramoto type population system in \cite{10}, the diffusion matrix is given by
\begin{equation}\label{sktdiffusion}
\begin{aligned}
A_{ij}(u)=\delta_{ij}p_i(u)+u_i\frac{\partial p_i}{\partial u_j}(u),\quad p_i(u)=a_{i0}+\sum_{k=1}^{n}a_{ik}u^s_k,\quad i,j=1,...,n.
\end{aligned}
\end{equation}

In (\ref{sktdiffusion}), $a_{i0},a_{ik}$ are positive constants, $s>0$ and $\delta_{ij}$ is the Kronecker delta symbol. If $s=1$, we choose the entropy as
\begin{equation}\label{sktentropy}
\begin{aligned}
h(u)=\sum_{i=1}^{n}\pi_iu_i(\log u_i-1),\quad\pi_i>0.
\end{aligned}
\end{equation}

Once diffusion coefficients $(a_{ij})$ satisfy certain assumptions, we can apply the entropy in (\ref{sktentropy}) to show a weak solution of the population system exists. In \cite{10}, assumptions for diffusion coefficients were discussed in detail. In general, we cannot show that $h^{\prime\prime}(u)A(u)$ in (\ref{sktdiffusion})-(\ref{sktentropy}) is positive semi-definite. This remains a great challenge in the study of parabolic partial differential equations. 

\section{Novelties and Main Steps}

We still rely on the entropy method to show the existence result. In this work, through the transformation of variables method, we show that a bounded weak solution of (\ref{equation})-(\ref{diffusiontwo}) exists.

The major novelty is we apply the transformation of variables technique. This allows us to derive a sequence of approximated solutions through the entropy method. We divide the proof into these steps.

\textbf{Step 1.} In Section 4, we first of all list (H1)-(H3) as assumptions for entropy density and diffusion matrix. These are major assumptions when we apply the entropy method in Section 5. 

\textbf{Step 2.} Based on these assumptions, we define $v=(v_1,v_2)$ in (\ref{v1v2}), $M(v)$ in (\ref{matrixv}) and $h(v)$ in (\ref{entropyepsilon}). We show that $h(v)$ satisfies (H1). Lemma \ref{positivesemidefinite} implies the positive definiteness of $h^{\prime\prime}(v)M(v)$.

\textbf{Step 3.} In Section 5, we prove Lemma \ref{mainresultv}, then to show the main result Theorem \ref{mainresult}. By the time discretization technique, we derive the approximated equations. Then we rely on (H1)-(H2) to derive a sequence of approximated solutions indexed by $\tau,\varepsilon>0$.

In this step, we mainly refer to Lemma 5, \cite{2}. For this classical lemma, we briefly state its idea of proof after we finish the proof of Theorem \ref{mainresult}.

\textbf{Step 4.} We provide uniform estimates of approximated solutions indexed by $\tau,\varepsilon>0$. By (H2)-(H3), we derive estimation formulas (\ref{firstaubin}) and (\ref{secondaubin}).

\textbf{Step 5.} We let $(\tau,\varepsilon)\to 0$. We apply the Aubin-Lions Lemma in \cite{17,18} to show that a weak solution of (\ref{equationv})-(\ref{matrixv}) exists. By the structural relationship between (\ref{equationv})-(\ref{matrixv}) and (\ref{equation})-(\ref{diffusiontwo}), we derive that a bounded weak solution of (\ref{equation})-(\ref{diffusiontwo}) exists.

\section{Preliminaries}

How to apply the entropy method to show weak solutions exist to cross-diffusion systems was discussed in \cite{2,7,8,9,10}. Let us present these assumptions:

(H1). There exists a convex function $h\in C^2(\mathcal{D},(0,\infty))$($\mathcal{D}\subset\mathbb{R}^n$ open and bounded, $n\geq 1$) such that its derivative $h^{\prime}:\mathcal{D}\to\mathbb{R}^n$ is invertible on $\mathbb{R}^n$;

(H2). For every $z=(z_1,...,z_n)^{\mathrm{T}}\in\mathbb{R}^n$ and $u=(u_1,...,u_n)\in\mathcal{D}$, there exists a constant $c>0$ that
\begin{equation}\nonumber
z^{\mathrm{T}}h^{\prime\prime}(u)A(u)z\geq c\sum_{i=1}^{n}z^2_i;
\end{equation}

(H3). There exists $a^{*}>0$ such that for all $u\in\mathcal{D}$ and $i,j=1,...,n$, it holds that $|A_{ij}(u)|<a^{*}$.

Let us define $v=(v_1,v_2)$, $\lambda=(\lambda_1,\lambda_2)$, such that 
\begin{equation}\label{v1v2}
\begin{aligned}
v_1=\lambda_1u_1,\quad v_2=\lambda_2u_2,\quad f(v)=\frac{\lambda_1\lambda_2}{\lambda_1\lambda_2+\mu_2\lambda_2v_1+\mu_1\lambda_1v_2},\quad\lambda_1,\lambda_2>0.
\end{aligned}
\end{equation}

We also denote 
\begin{equation}\nonumber
\begin{aligned}
&\mathcal{P}_{\lambda}=\{u=(u_1,u_2):\lambda_1u_1+\lambda_2u_2<1,u_1,u_2>0\},\\&\mathcal{D}=\{v=(v_1,v_2):v_1+v_2<1,v_1,v_2>0\}.
\end{aligned}
\end{equation}
By these definitions, $u\in\mathcal{P}_{\lambda}$ if and only if $v\in\mathcal{D}$. 

In the following discussion, we consider equations
\begin{equation}\label{equationv}
\partial_tv=\dive\big(M(v)\nabla v\big)\text{ in }Q_T=\Omega\times(0,T),\text{ }T>0,
\end{equation}
and
\begin{equation}\label{initialv}
\begin{aligned}
&(M_{11}(v)\nabla v_1+M_{12}(v)\nabla v_2)\cdot\nu=0,\quad(M_{21}(v)\nabla v_1+M_{22}(v)\nabla v_2)\cdot\nu=0\text{ on }\partial\Omega,\\&v^0_1=\lambda_1u^0_1,\quad v^0_2=\lambda_2u^0_2\text{ in }\Omega,
\end{aligned}
\end{equation}
with the diffusion matrix $M(v)=(M_{ij}(v))$ given by
\begin{equation}\label{matrixv}
\begin{aligned}
&M(v)/f(v)=\left(\begin{array}{cccccc}
\mu_1 & 0 \\
0 & \mu_2 \\
\end{array}\right)+\mu_1\mu_2\left(\begin{array}{cccccc}
v_1/\lambda_1 & v_1/\lambda_2 \\
v_2/\lambda_1 & v_2/\lambda_2 \\
\end{array}\right).
\end{aligned}
\end{equation}

Let us give the definition of weak solutions to (\ref{equationv})-(\ref{matrixv}).
\begin{definition}\label{defweakv}
For every $T>0$, $v^0=(v^0_1,v^0_2)\in\mathcal{D}$, if $v=(v_1,v_2)$ satisfies that $v\in\overline{\mathcal{D}}$, and for all test functions $\phi=(\phi_1,\phi_2)$, $\phi_1,\phi_2\in L^2(0,T;H^1(\Omega))$,
\begin{equation}\nonumber
\begin{aligned}
&\int_{0}^{T}\langle\partial_tv_1,\phi_1\rangle dt+\int_{Q_T}\frac{\mu_1\lambda_2(\lambda_1+\mu_2v_1)\nabla v_1+\mu_1\mu_2\lambda_1v_1\nabla v_2}{\lambda_1\lambda_2+\mu_2\lambda_2v_1+\mu_1\lambda_1v_2}\cdot\nabla\phi_1dxdt=0,\\&\int_{0}^{T}\langle\partial_tv_2,\phi_2\rangle dt+\int_{Q_T}\frac{\mu_1\mu_2\lambda_2v_2\nabla v_1+\mu_2\lambda_1(\lambda_2+\mu_1v_2)\nabla v_2}{\lambda_1\lambda_2+\mu_2\lambda_2v_1+\mu_1\lambda_1v_2}\cdot\nabla\phi_2dxdt=0,
\end{aligned}
\end{equation}
with the property that $v_i\in L^2(0,T;H^1(\Omega))$, $\partial_tv_i\in L^2(0,T;H^1(\Omega)^{\prime})$, $i=1,2$, then we call $v$ is a weak solution of $(\ref{equationv})$-$(\ref{matrixv})$.
\end{definition}

In Lemma \ref{uvtransfer}, we consider the relationship between $u=(u_1,u_2)$ and $v=(v_1,v_2)$.

\begin{lem}\label{uvtransfer}
If $v^0=(v^0_1,v^0_2)\in\mathcal{D}$, $v=(v_1,v_2)\in\overline{\mathcal{D}}$ is a weak solution of $(\ref{equationv})$-$(\ref{matrixv})$, then $u^0=(u^0_1,u^0_2)=(v^0_1/\lambda_1,v^0_2/\lambda_2)\in\mathcal{P}_{\lambda}$, $u=(u_1,u_2)\in\overline{\mathcal{P}}_{\lambda}$. For all test functions $\phi=(\phi_1,\phi_2)$, $\phi_1,\phi_2\in L^2(0,T;H^1(\Omega))$, $u=(u_1,u_2)$ satisfies that
\begin{equation}\label{usatisfy}
\begin{aligned}
&\int_{0}^{T}\langle\partial_tu_1,\phi_1\rangle dt+\int_{Q_T}\frac{\mu_1(1+\mu_2u_1)\nabla u_1+\mu_1\mu_2u_1\nabla u_2}{1+\mu_2u_1+\mu_1u_2}\cdot\nabla\phi_1dxdt=0,\\&\int_{0}^{T}\langle\partial_tu_2,\phi_2\rangle dt+\int_{Q_T}\frac{\mu_2(1+\mu_1u_2)\nabla u_2+\mu_1\mu_2u_2\nabla u_1}{1+\mu_2u_1+\mu_1u_2}\cdot\nabla\phi_2dxdt=0,
\end{aligned}
\end{equation}
with the property that $u_i\in L^2(0,T;H^1(\Omega))$, $\partial_tu_i\in L^2(0,T;H^1(\Omega)^{\prime})$, $i=1,2$.
\end{lem}

\begin{proof}

Equations in (\ref{equationv})-(\ref{matrixv}) are
\begin{equation}\nonumber
\begin{aligned}
&\partial_tv_1=\dive\Big(f(v)\cdot\big(\mu_1\nabla v_1+\mu_1\mu_2(\frac{v_1}{\lambda_1}\nabla v_1+\frac{v_1}{\lambda_2}\nabla v_2)\big)\Big),
\\&\partial_tv_2=\dive\Big(f(v)\cdot\big(\mu_2\nabla v_2+\mu_1\mu_2(\frac{v_2}{\lambda_1}\nabla v_1+\frac{v_2}{\lambda_2}\nabla v_2)\big)\Big),
\end{aligned}
\end{equation}
which are equivalent to
\begin{equation}\label{equivalentv}
\begin{aligned}
&\partial_t\big(\frac{v_1}{\lambda_1}\big)=\dive\Big(f(v)\cdot\Big(\mu_1\nabla\big(\frac{v_1}{\lambda_1}\big)+\mu_1\mu_2\big(\frac{v_1}{\lambda_1}\nabla\big(\frac{v_1}{\lambda_1}\big)+\frac{v_1}{\lambda_1}\nabla\big(\frac{v_2}{\lambda_2}\big)\big)\Big)\Big),
\\&\partial_t\big(\frac{v_2}{\lambda_2}\big)=\dive\Big(f(v)\cdot\Big(\mu_2\nabla\big(\frac{v_2}{\lambda_2}\big)+\mu_1\mu_2\big(\frac{v_2}{\lambda_2}\nabla\big(\frac{v_1}{\lambda_1}\big)+\frac{v_2}{\lambda_2}\nabla\big(\frac{v_2}{\lambda_2}\big)\big)\Big)\Big).
\end{aligned}
\end{equation}

By (\ref{v1v2}), we conclude that (\ref{equivalentv}) is
\begin{equation}\label{equivalentu}
\begin{aligned}
&\partial_tu_1=\dive\Big(\frac{1}{1+\mu_2u_1+\mu_1u_2}\big(\mu_1\nabla u_1+\mu_1\mu_2(u_1\nabla u_1+u_1\nabla u_2)\big)\Big),
\\&\partial_tu_2=\dive\Big(\frac{1}{1+\mu_2u_1+\mu_1u_2}\big(\mu_2\nabla u_2+\mu_1\mu_2(u_2\nabla u_1+u_2\nabla u_2)\big)\Big).
\end{aligned}
\end{equation}

By (\ref{equivalentu}), we deduce that once $v=(v_1,v_2)$ is a weak solution of (\ref{equationv})-(\ref{matrixv}), then $u^0=(u^0_1,u^0_2)\in\mathcal{P}_{\lambda}$, $u=(u_1,u_2)\in\overline{\mathcal{P}}_{\lambda}$, and satisfy (\ref{usatisfy}). 

\end{proof}

Let us denote $v_3=1-v_1-v_2$, and define the entropy $h(v)$ as
\begin{equation}\label{entropyepsilon}
\begin{aligned}
&h(v)=v_1(\log v_1-1)+v_2(\log v_2-1)+v_3(\log v_3-1)+3.
\end{aligned}
\end{equation}

We notice that $h\in C^2(\mathcal{D},(0,\infty))$, and
\begin{equation}\nonumber
\begin{aligned}
\partial h/\partial v_i=\log v_i-\log v_3,\quad i=1,2.
\end{aligned}
\end{equation}

Let us define 
\begin{equation}\nonumber
\begin{aligned}
g(y)=(e^{w_1}+e^{w_2})\cdot(1-y),\quad 0<y<1.
\end{aligned}
\end{equation}
$g$ is continuous and nonincreasing in $0<y<1$, and for every $w=(w_1,w_2)\in\mathbb{R}^2$,
\begin{equation}\nonumber
\begin{aligned}
g(0)=e^{w_1}+e^{w_2}>0,\quad g(1)=0.
\end{aligned}
\end{equation}
We deduce that there exists a unique fixed point $0<y_0<1$, such that $g(y_0)=y_0$. 

Then we define $v_i=e^{w_i}(1-y_0)>0$($i=1,2$), $0<y_0<1$, which satisfies $v_1+v_2=g(y_0)=y_0<1$. Consequently, $v=(v_1,v_2)\in\mathcal{D}$, $v_3=1-v_1-v_2=1-y_0$, and $$w_i=\log v_i-\log v_3=\partial h/\partial v_i,\quad i=1,2.$$ Then we show that $h^{\prime}:\mathcal{D}\to\mathbb{R}^2$ is invertible on $\mathbb{R}^2$, and assumption (H1) can be satisfied.

Let us denote $$c=\frac{\lambda_1\lambda_2\min\{\mu_1,\mu_2\}}{\lambda_1\lambda_2+\max\{\mu_1,\mu_2\}(\lambda_1+\lambda_2)}.$$ In Lemma \ref{positivesemidefinite}, we verify assumption (H2).

\begin{lem}\label{positivesemidefinite}
If constants $\lambda_1,\lambda_2>0$ satisfy the relation
\begin{equation}\label{lambdaonetwo}
\begin{aligned}
\mu_1+\mu_1\mu_2/\lambda_1=\mu_2+\mu_1\mu_2/\lambda_2,
\end{aligned}
\end{equation}
then for every $z=(z_1,z_2)^{\mathrm{T}}\in\mathbb{R}^2$ and $v\in\mathcal{D}$, we have
\begin{equation}\nonumber
\begin{aligned}
&z^{\mathrm{T}}h^{\prime\prime}(v)M(v)z\geq c(z^2_1+z^2_2).
\end{aligned}
\end{equation}
\end{lem}

\begin{proof}

By (\ref{entropyepsilon}), we have for every $v\in\mathcal{D}$,
\begin{equation}\label{hessianmatrix}
\begin{aligned}
&h^{\prime\prime}(v)=\left(\begin{array}{cccccc}
\frac{1}{v_1} & 0 \\
 0 & \frac{1}{v_2} \\
\end{array}\right)+\frac{1}{v_3}\left(\begin{array}{cccccc}
1 & 1 \\
1 & 1 \\
\end{array}\right).
\end{aligned}
\end{equation}

We multiply $M(v)$ by $h^{\prime\prime}(v)$, which is
\begin{equation}\nonumber
\begin{aligned}
&h^{\prime\prime}(v)M(v)/f(v)=\Big\{\left(\begin{array}{cccccc}
\frac{1}{v_1} & 0 \\
 0 & \frac{1}{v_2} \\
\end{array}\right)+\frac{1}{v_3}\left(\begin{array}{cccccc}
1 & 1 \\
1 & 1 \\
\end{array}\right)\Big\}\cdot\Big\{\left(\begin{array}{cccccc}
\mu_1 & 0 \\
0 & \mu_2 \\
\end{array}\right)+\mu_1\mu_2\left(\begin{array}{cccccc}
v_1/\lambda_1 & v_1/\lambda_2 \\
v_2/\lambda_1 & v_2/\lambda_2 \\
\end{array}\right)\Big\}\\&=\left(\begin{array}{cccccc}
\frac{\mu_1}{v_1} & 0 \\
0 & \frac{\mu_2}{v_2} \\
\end{array}\right)+\frac{1}{v_3}\left(\begin{array}{cccccc}
\mu_1 & \mu_2 \\
\mu_1 & \mu_2\\
\end{array}\right)+\mu_1\mu_2\left(\begin{array}{cccccc}
1/\lambda_1 & 1/\lambda_2 \\
1/\lambda_1 & 1/\lambda_2 \\
\end{array}\right)+\frac{\mu_1\mu_2(v_1+v_2)}{v_3}\left(\begin{array}{cccccc}
1/\lambda_1 & 1/\lambda_2 \\
1/\lambda_1 & 1/\lambda_2 \\
\end{array}\right).
\end{aligned}
\end{equation}

By the relation (\ref{lambdaonetwo}), we conclude that
\begin{equation}\nonumber
\begin{aligned}
h^{\prime\prime}(v)M(v)/f(v)&=\left(\begin{array}{cccccc}
\frac{\mu_1}{v_1} & 0 \\
0 & \frac{\mu_2}{v_2} \\
\end{array}\right)+\frac{1}{v_3}\left(\begin{array}{cccccc}
\mu_1+\mu_1\mu_2/\lambda_1 & \mu_2+\mu_1\mu_2/\lambda_2 \\
\mu_1+\mu_1\mu_2/\lambda_1 & \mu_2+\mu_1\mu_2/\lambda_2 \\
\end{array}\right)\\&=\left(\begin{array}{cccccc}
\frac{\mu_1}{v_1} & 0 \\
0 & \frac{\mu_2}{v_2} \\
\end{array}\right)+\frac{\mu_1+\mu_1\mu_2/\lambda_1}{v_3}\left(\begin{array}{cccccc}
1 & 1 \\
1 & 1 \\
\end{array}\right).
\end{aligned}
\end{equation}

For $v\in\mathcal{D}$, $z=(z_1,z_2)^{\mathrm{T}}\in\mathbb{R}^2$, we have
\begin{equation}\nonumber
\begin{aligned}
&z^{\mathrm{T}}h^{\prime\prime}(v)M(v)z=f(v)\cdot\Big(\frac{\mu_1z^2_1}{v_1}+\frac{\mu_2z^2_2}{v_2}+\frac{\mu_1+\mu_1\mu_2/\lambda_1}{v_3}\cdot(z_1+z_2)^2\Big),
\end{aligned}
\end{equation}
which follows that
\begin{equation}\nonumber
\begin{aligned}
&z^{\mathrm{T}}h^{\prime\prime}(v)M(v)z\geq c(z^2_1+z^2_2),
\end{aligned}
\end{equation}
and we finish the proof of Lemma \ref{positivesemidefinite}.
\end{proof}

By Lemma \ref{positivesemidefinite}, $h(v)$ and $A(v)$ satisfy (H2). In Lemma \ref{lambdam}, we discuss the structural relationship between $\mathcal{D}_M$ and $\mathcal{P}_{\lambda}$.

\begin{lem}\label{lambdam}
For every $M>0$, there corresponds $\lambda_1,\lambda_2>0$ satisfying the relation $(\ref{lambdaonetwo})$, and if $u^0=(u^0_1,u^0_2)\in\mathcal{D}_{M}$, then $u^0\in\mathcal{P}_{\lambda}$.
\end{lem}

\begin{proof}
Once constants $\lambda_1,\lambda_2$ satisfy (\ref{lambdaonetwo}), and $\mu_1>\mu_2$, then we have $$\lambda_1\lambda_2(\mu_1-\mu_2)=\mu_1\mu_2(\lambda_1-\lambda_2),$$ which follows 
\begin{equation}\label{mu1mu2}
\begin{aligned}
\lambda_2=\frac{\mu_1\mu_2\lambda_1}{(\mu_1-\mu_2)\lambda_1+\mu_1\mu_2}.
\end{aligned}
\end{equation}

Let us choose $0<\lambda_1<1/M$, then by (\ref{mu1mu2}), we have $0<\lambda_2<1/M$.

Once constants $\lambda_1,\lambda_2$ satisfy (\ref{lambdaonetwo}), and $\mu_1<\mu_2$, then we have $$\lambda_1\lambda_2(\mu_2-\mu_1)=\mu_1\mu_2(\lambda_2-\lambda_1),$$ which follows 
\begin{equation}\label{mu2mu1}
\begin{aligned}
\lambda_1=\frac{\mu_1\mu_2\lambda_2}{(\mu_2-\mu_1)\lambda_2+\mu_1\mu_2}.
\end{aligned}
\end{equation}

Let us choose $0<\lambda_2<1/M$, then by (\ref{mu2mu1}), we have $0<\lambda_1<1/M$. We conclude that we are allowed to find $\lambda_1,\lambda_2$ satisfying (\ref{lambdaonetwo}) and $0<\lambda_1,\lambda_2<1/M$ simultaneously.

If $u^0\in\mathcal{D}_M$, then $u^0_1,u^0_2>0$, and
\begin{equation}\label{mu2mu10}
\begin{aligned}
\lambda_1u^0_1+\lambda_2u^0_2<u^0_1/M+u^0_2/M<1.
\end{aligned}
\end{equation}

By (\ref{mu2mu10}), we deduce that once $u^0\in\mathcal{D}_{M}$, then $u^0\in\mathcal{P}_{\lambda}$, and we finish the proof.

\end{proof}

These preliminary results are key factors in the proof of Lemma \ref{mainresultv} and Theorem \ref{mainresult}.

\section{Proof of the Main Theorem}

Let us present Lemma \ref{mainresultv} and give the proof.
\begin{lem}\label{mainresultv}
Let $v^0=(v^0_1,v^0_2)$ be such that $v^0\in\mathcal{D}$. If constants $\lambda_1,\lambda_2>0$ satisfy relation $(\ref{lambdaonetwo})$, then there exists a weak solution $v\in\overline{\mathcal{D}}$ of $(\ref{equationv})$-$(\ref{matrixv})$.
\end{lem}

\begin{proof}

\textbf{Step 1.} Let $T>0$, $N\in\mathbb{N}$, $\tau=T/N$, $0<\varepsilon<1$, and $m\in\mathbb{N}$ with $m>d/2$. This ensures that the embedding $H^m(\Omega)\hookrightarrow L^{\infty}(\Omega)$ is compact. 

Lemma \ref{positivesemidefinite} indicates that assumption (H2) can be satisfied. In \cite{2}, Lemma 5 showed that given $w^{k-1}\in L^{\infty}(\Omega)$ for $k\in\mathbb{N}$, there exists $w^k\in H^m(\Omega)$, $v(w^k(x))\in\mathcal{D}$ for $x\in\Omega$ that
\begin{equation}\label{inductionhyp}
\begin{aligned}
&\frac{1}{\tau}\int_{\Omega}(v(w^k)-v(w^{k-1}))\cdot\phi dx+\int_{\Omega}\nabla\phi:B(w^k)\nabla w^kdx\\&+\varepsilon\int_{\Omega}\Big(\sum_{|\alpha|=m}D^{\alpha}w^k\cdot D^{\alpha}\phi+w^k\cdot\phi\Big)dx=0,
\end{aligned}
\end{equation}
for every $\phi\in H^m(\Omega)$. Here, $A:B=\sum_{i,j}a_{ij}b_{ij}$ for matrices $A=(a_{ij})$ and $B=(b_{ij})$.

In (\ref{inductionhyp}), $v(w^k)=(h^{\prime})^{-1}(w^k)$, $B(w^k)=M(v(w^k))H(v(w^k))^{-1}$, and $H(v)=h^{\prime\prime}(v)$. $|\alpha|=\alpha_1+\cdot\cdot\cdot+\alpha_d=m$ is a multiindex and $D^{\alpha}=\partial^{|\alpha|}/(\partial x^{\alpha_1}_1\cdot\cdot\cdot\partial x^{\alpha_d}_d)$ is a partial derivative of order $m$. If $k=1$, we define $w^0=h^{\prime}(v^0)$.

Let us denote $v^k=v(w^k)$, and introduce the piecewise in time constant functions $v^{(\tau)}=v^k$, $w^{(\tau)}=w^k$, for $x\in\Omega$, $t\in((k-1)\tau,k\tau]$. At time $t=0$, we set $v^{(\tau)}(\cdot,0)=(h^{\prime})^{-1}(v^0)$. Let $v^{(\tau)}=(v^{(\tau)}_1,v^{(\tau)}_2)$, we define the backward shift operator $(\sigma_{\tau}v^{(\tau)})(x,t)=v^{k-1}$ for $x\in\Omega$, $t\in((k-1)\tau,k\tau]$. $v^{(\tau)}$ satisfies that
\begin{equation}\label{inductiontau}
\begin{aligned}
&\frac{1}{\tau}\int_{0}^{T}\int_{\Omega}(v^{(\tau)}-\sigma_{\tau}v^{(\tau)})\cdot\phi dxdt+\int_{0}^{T}\int_{\Omega}\nabla\phi:B(w^{(\tau)})\nabla w^{(\tau)}dxdt\\&+\varepsilon\int_{0}^{T}\int_{\Omega}\Big(\sum_{|\alpha|=m}D^{\alpha}w^{(\tau)}\cdot D^{\alpha}\phi+w^{(\tau)}\cdot\phi\Big)dxdt=0,
\end{aligned}
\end{equation}
for piecewise constant functions $\phi:(0,T)\to H^m(\Omega)$. By a density argument, the equation (\ref{inductiontau}) also holds for every $\phi\in L^2(0,T;H^m(\Omega))$.

Lemma 5, \cite{2} only requires that $h^{\prime\prime}(v)M(v)$ is positive semi-definite. The positive semi-definiteness of $h^{\prime\prime}(v)M(v)$ is included in (H2). 

Lemma 5, \cite{2} showed the existence of approximated solutions by applying Lax-Milgram lemma and fixed point theorem. By Lemma 5 of \cite{2}, we derive a sequence of approximated solutions and the entropy inequality for these solutions.

The proof of Lemma 5, \cite{2} takes several pages and this classical result has been widely applied in global existence analysis of cross-diffusion systems. For these reasons, we briefly state the proof of this classical Lemma 5 after we finish the proof of Theorem \ref{mainresult}.

\textbf{Step 2.} Relying on (\ref{inductionhyp}) and the inequality
\begin{equation}\nonumber
\begin{aligned}
h(v(w^k))-h(v(w^{k-1}))\leq h^{\prime}(v(w^k))\cdot(v(w^k)-v(w^{k-1})),
\end{aligned}
\end{equation}
we deduce that $w^k$ satisfies the discrete entropy inequality
\begin{equation}\label{entropyinequality}
\begin{aligned}
&\int_{\Omega}h(v(w^k))dx+\tau\int_{\Omega}\nabla w^k:B(w^k)\nabla w^kdx\\&+\varepsilon\tau\int_{\Omega}\Big(\sum_{|\alpha|=m}|D^{\alpha}w^k|^2+|w^k|^2\Big)dx\leq\int_{\Omega}h(v(w^{k-1}))dx.
\end{aligned}
\end{equation}

By the generalized Poincar\'e inequality, we deduce that
\begin{equation}\nonumber
\begin{aligned}
\int_{\Omega}\Big(\sum_{|\alpha|=m}|D^{\alpha}w^k|^2+|w^k|^2\Big)dx\geq C_p\|w^k\|^2_{H^m(\Omega)},
\end{aligned}
\end{equation}
where $C_p>0$ is the Poincar\'e constant. By Lemma \ref{positivesemidefinite}, we have
\begin{equation}\nonumber
\begin{aligned}
&\nabla w^k:B(w^k)\nabla w^k=\nabla v^k:H(v^k)M(v^k)\nabla v^k\geq c(|\nabla v^{k}_{1}|^2+|\nabla v^{k}_2|^2).
\end{aligned}
\end{equation}

Combining above computation results, by (\ref{entropyinequality}), we conclude that
\begin{equation}\label{konetojw}
\begin{aligned}
&\int_{\Omega}h(v^k)dx+c\tau\int_{\Omega}|\nabla v^k|^2dx+\varepsilon C_p\tau\|w^k\|^2_{H^m(\Omega)}\leq\int_{\Omega}h(v^{k-1})dx.
\end{aligned}
\end{equation}

Summing these inequalities over $k=1,...,j$ in (\ref{konetojw}), we can show that
\begin{equation}\nonumber
\begin{aligned}
&\int_{\Omega}h(v^j)dx+c\tau\sum_{k=1}^{j}\int_{\Omega}|\nabla v^k|^2dx+\varepsilon C_p\tau\sum_{k=1}^{j}\|w^k\|^2_{H^m(\Omega)}\leq\int_{\Omega}h(v^0)dx,
\end{aligned}
\end{equation}
then we derive that
\begin{equation}\label{firstaubin}
\begin{aligned}
&\|v^{(\tau)}\|_{L^2(0,T;H^1(\Omega))}+\varepsilon^{1/2}\|w^{(\tau)}\|_{L^2(0,T;H^m(\Omega))}\leq C,
\end{aligned}
\end{equation}
with $C>0$ a generic constant independent of $\tau,\varepsilon$. 

\textbf{Step 3.} Since $M(v)$ is uniformly bounded in $v\in\mathcal{D}$, then by (\ref{firstaubin}), we can show that
\begin{equation}\label{abounded}
\begin{aligned}
\Big|\int_{0}^{T}\int_{\Omega}\nabla\phi:M(v^{(\tau)})\nabla v^{(\tau)}dxdt\Big|\leq C.
\end{aligned}
\end{equation}

In (\ref{inductiontau}), we have
\begin{equation}\label{secondaubin}
\begin{aligned}
&\Big|\frac{1}{\tau}\int_{0}^{T}\int_{\Omega}(v^{(\tau)}-\sigma_{\tau}v^{(\tau)})\cdot\phi dxdt\Big|\leq\Big|\int_{0}^{T}\int_{\Omega}\nabla\phi:M(v^{(\tau)})\nabla v^{(\tau)}dxdt\Big|\\&+\Big|\varepsilon\int_{0}^{T}\int_{\Omega}\Big(\sum_{|\alpha|=m}D^{\alpha}w^{(\tau)}\cdot D^{\alpha}\phi+w^{(\tau)}\cdot\phi\Big)dxdt\Big|\leq C.
\end{aligned}
\end{equation}

By (\ref{firstaubin}) and (\ref{secondaubin}), we apply the Aubin-Lions Lemma. There exists a subsequence which we do not relabel, such that as $(\tau,\varepsilon)\to 0$,
\begin{equation}\label{firstconvergence}
\begin{aligned}
v^{(\tau)}&\to v\text{ strongly in }L^2(0,T;L^2(\Omega)),\text{ and a.e. in }\Omega\times(0,T).
\\v^{(\tau)}&\rightharpoonup v\text{ weakly in }L^2(0,T;H^1(\Omega)).
\\\tau^{-1}(v^{(\tau)}-\sigma_{\tau}v^{(\tau)})&\rightharpoonup\partial_tv\text{ weakly in }L^2(0,T;H^m(\Omega)^{\prime}).\\\varepsilon w^{(\tau)}&\to 0\text{ strongly in }L^2(0,T;H^m(\Omega)).
\end{aligned}
\end{equation}

Let us denote
\begin{equation}\nonumber
\begin{aligned}
&G_1=f(v^{(\tau)})(\mu_1\nabla v^{(\tau)}_1+(\mu_1\mu_2/\lambda_1)v^{(\tau)}_1\nabla v^{(\tau)}_1+(\mu_1\mu_2/\lambda_2)v^{(\tau)}_1\nabla v^{(\tau)}_2),
\\&G_2=f(v^{(\tau)})(\mu_2\nabla v^{(\tau)}_2+(\mu_1\mu_2/\lambda_1)v^{(\tau)}_2\nabla v^{(\tau)}_1+(\mu_1\mu_2/\lambda_2)v^{(\tau)}_2\nabla v^{(\tau)}_2),
\\&I_1=f(v)(\mu_1\nabla v_1+(\mu_1\mu_2/\lambda_1)v_1\nabla v_1+(\mu_1\mu_2/\lambda_2)v_1\nabla v_2),
\\&I_2=f(v)(\mu_2\nabla v_2+(\mu_1\mu_2/\lambda_1)v_2\nabla v_1+(\mu_1\mu_2/\lambda_2)v_2\nabla v_2),
\end{aligned}
\end{equation}
and $J=(J_1,J_2)$, with
\begin{equation}\nonumber
\begin{aligned}
&J_1=f(v^{(\tau)})\big(\mu_1\nabla(v^{(\tau)}_1-v_1)+(\mu_1\mu_2/\lambda_1)(v^{(\tau)}_1\nabla v^{(\tau)}_1-v_1\nabla v_1)+(\mu_1\mu_2/\lambda_2)(v^{(\tau)}_1\nabla v^{(\tau)}_2-v_1\nabla v_2)\big)\\&+(f(v^{(\tau)})-f(v))\cdot(\mu_1\nabla v_1+(\mu_1\mu_2/\lambda_1)v_1\nabla v_1+(\mu_1\mu_2/\lambda_2)v_1\nabla v_2),
\\&J_2=f(v^{(\tau)})\big(\mu_2\nabla(v^{(\tau)}_2-v_2)+(\mu_1\mu_2/\lambda_1)(v^{(\tau)}_2\nabla v^{(\tau)}_1-v_2\nabla v_1)+(\mu_1\mu_2/\lambda_2)(v^{(\tau)}_2\nabla v^{(\tau)}_2-v_2\nabla v_2)\big)\\&+(f(v^{(\tau)})-f(v))\cdot(\mu_2\nabla v_2+(\mu_1\mu_2/\lambda_1)v_2\nabla v_1+(\mu_1\mu_2/\lambda_2)v_2\nabla v_2).
\end{aligned}
\end{equation}

We have $G_1-I_1=J_1$, $G_2-I_2=J_2$. For a.e. $(x,t)\in\Omega\times(0,T)$, as $(\tau,\varepsilon)\to 0$,
\begin{equation}\label{aeconvergence}
\begin{aligned}
v^{(\tau)}_1\to v_1,\quad v^{(\tau)}_2\to v_2,\quad f(v^{(\tau)})\to f(v).
\end{aligned}
\end{equation}

We notice that 
\begin{equation}\label{u1u2u3difference}
\begin{aligned}
v^{(\tau)}_i\nabla v^{(\tau)}_j-v_i\nabla v_j=v^{(\tau)}_i\nabla(v^{(\tau)}_j-v_j)+(v^{(\tau)}_i-v_i)\nabla v_j,\quad 1\leq i,j\leq 2.
\end{aligned}
\end{equation}
Then by (\ref{firstconvergence})-(\ref{u1u2u3difference}) and dominated convergence theorem, we can show that for every $\phi\in L^2(0,T;H^m(\Omega))$, as $(\tau,\varepsilon)\to 0$,
\begin{equation}\label{j1j3}
\begin{aligned}
\int_{0}^{T}\int_{\Omega}J\cdot\nabla\phi dxdt\to 0.
\end{aligned}
\end{equation}

By (\ref{j1j3}), there exists a subsequence which we do not relabel, such that for every $\phi\in L^2(0,T;H^m(\Omega))$, as $(\tau,\varepsilon)\to 0$, 
\begin{equation}\label{secondconvergence}
\begin{aligned}
&\int_{0}^{T}\int_{\Omega}\big(M(v^{(\tau)})\nabla v^{(\tau)}-M(v)\nabla v\big)\cdot\nabla\phi dxdt\to 0.
\end{aligned}
\end{equation}

By (\ref{firstconvergence}) and (\ref{secondconvergence}), we conclude that
\begin{equation}\label{csurplus}
\begin{aligned}
\int_{0}^{T}\langle\partial_tv,\phi\rangle dt+\int_{0}^{T}\int_{\Omega}\nabla\phi:M(v)\nabla vdxdt=0.
\end{aligned}
\end{equation}

By (\ref{secondaubin}), $\partial_tv\in L^2(0,T;H^1(\Omega)^{\prime})$, and consequently, the weak formulation (\ref{csurplus}) also holds for every $\phi\in L^2(0,T;H^1(\Omega))$. Finally, $v\in H^1(0,T;H^1(\Omega)^{\prime})\hookrightarrow C^0([0,T];H^1(\Omega)^{\prime})$ shows that the initial condition is satisfied in $H^1(\Omega)^{\prime}$. $v$ satisfies (\ref{equationv})-(\ref{matrixv}) in a weak sense, and we finish the proof.

\end{proof}

After all these preparations, let us give the proof of the main result Theorem \ref{mainresult}.

\begin{proof}

By Lemma \ref{lambdam}, we choose $\lambda_1,\lambda_2>0$ satisfying the relation (\ref{lambdaonetwo}), and $\mathcal{D}_{M}\subset\mathcal{P}_{\lambda}$. By Lemma \ref{mainresultv}, a weak solution $v=(v_1,v_2)$ of (\ref{equationv})-(\ref{matrixv}) exists. 

By Lemma \ref{uvtransfer}, since $u^0\in\mathcal{D}_M\subset\mathcal{P}_{\lambda}$, then a weak solution $u=(u_1,u_2)\in\overline{\mathcal{P}}_{\lambda}\subset\overline{\mathcal{D}}_{M^{\prime}}$ of (\ref{equation})-(\ref{diffusiontwo}) exists. In here, $M^{\prime}>M$ is a constant. Then we have shown Theorem \ref{mainresult}.

\end{proof}

For integrity purpose, we briefly state Lemma 5 of \cite{2} in the end of this section. Let $y\in L^{\infty}(\Omega;\mathbb{R}^n)$ and $\delta\in[0,1]$,  we consider the equation
\begin{equation}\label{laxmilgram}
a(w,\phi)=F(\phi),\text{ for all }\phi\in H^m(\Omega;\mathbb{R}^n),
\end{equation}
where
\begin{equation}\nonumber
\begin{aligned}
&a(w,\phi)=\int_{\Omega}\nabla\phi:B(y)\nabla wdx+\varepsilon\int_{\Omega}\Big(\sum_{|\alpha|=m}D^{\alpha}w\cdot D^{\alpha}\phi+w\cdot\phi\Big)dx,\\&F(\phi)=-\frac{\delta}{\tau}\int_{\Omega}(v(y)-v(w^{k-1}))\cdot\phi dx.
\end{aligned}
\end{equation}

The forms $a$ and $F$ are bounded on $H^m(\Omega;\mathbb{R}^n)$. The matrix $B(y)=M(v(y))(h^{\prime\prime}(v(y)))^{-1}$ is positive semi-definite since
\begin{equation}\nonumber
\begin{aligned}
z^{\mathrm{T}}B(y)z=((h^{\prime\prime}(v))^{-1}z)^{\mathrm{T}}h^{\prime\prime}(v)M(v(y))(h^{\prime\prime}(v))^{-1}z\geq 0,\text{ for all }z\in\mathbb{R}^n.
\end{aligned}
\end{equation}

By the Poincar\'e inequality, the bilinear form $a$ is coercive:
\begin{equation}\nonumber
\begin{aligned}
a(w,w)\geq\varepsilon\int_{\Omega}\Big(\sum_{|\alpha|=m}|D^{\alpha}w|^2+|w|^2\Big)dx\geq\varepsilon C_p\|w\|^2_{H^m(\Omega)}.
\end{aligned}
\end{equation}

We apply the Lax-Milgram lemma to obtain the existence of a unique solution $w\in H^m(\Omega;\mathbb{R}^n)\hookrightarrow L^{\infty}(\Omega;\mathbb{R}^n)$ to (\ref{laxmilgram}). This defines the fixed point operator $S:L^{\infty}(\Omega;\mathbb{R}^n)\times[0,1]\to L^{\infty}(\Omega;\mathbb{R}^n)$, $S(y,\delta)=w$, where $w$ solves (\ref{laxmilgram}).

The fixed point operator $S$ is continuous, compact and bounded. For continuity and compactness of the operator $S$, please refer to \cite{2} for a detailed explanation. It remains to prove a uniform bound for all fixed points of $S(\cdot,\delta)$ in $L^{\infty}(\Omega;\mathbb{R}^n)$. Let $w\in L^{\infty}(\Omega;\mathbb{R}^n)$ be such a fixed point. Then $w$ solves (\ref{laxmilgram}) with $y$ replaced by $w$. With the test function $\phi=w$, we find that
\begin{equation}\nonumber
\begin{aligned}
\frac{\delta}{\tau}\int_{\Omega}(v(w)-v(w^{k-1}))\cdot wdx+\int_{\Omega}\nabla w:B(w)\nabla wdx+\varepsilon\int_{\Omega}\Big(\sum_{|\alpha|=m}|D^{\alpha}w|^2+|w|^2\Big)dx=0.
\end{aligned}
\end{equation}

The convexity of $h$ implies that $h(x)-h(y)\leq h^{\prime}(x)(x-y)$ for all $x,y\in\mathcal{D}$. Let us choose $x=v(w),y=v(w^{k-1})$ and $h^{\prime}(v(w))=w$, this gives
\begin{equation}\nonumber
\begin{aligned}
\delta\int_{\Omega}h(v(w))dx+\varepsilon\tau\int_{\Omega}\Big(\sum_{|\alpha|=m}|D^{\alpha}w|^2+|w|^2\Big)dx\leq\delta\int_{\Omega}h(v(w^{k-1}))dx.
\end{aligned}
\end{equation}

This yields an $H^m$ bound for $w$ uniform in $\delta$. The Leray-Schauder fixed point theorem shows that there exists a solution $w\in H^m(\Omega;\mathbb{R}^n)$ to (\ref{laxmilgram}), with $y$ replaced by $w$ and $\delta=1$.

\section{Conclusions}

The entropy method is a powerful tool in global existence analysis of cross-diffusion systems. Once the cross-diffusion system possesses the entropy structure in \cite{19}, a systematic proof scheme can be applied to show a weak solution exists. How to show weak solutions exist to cross-diffusion systems without entropy structures is a central topic in the study of parabolic partial differential equations.

In semiconductor model (\ref{equation})-(\ref{diffusiontwo}), if $\mu_1\neq\mu_2$, $h^{\prime\prime}_2(u)A(u)$ is not positive semi-definite. We apply the transformation of variables technique to show the main result Theorem \ref{mainresult}. The conclusion is semiconductor model has a bounded weak solution. In \cite{20}, bounded weak solutions to a class of cross-diffusion systems were also discussed. In the proof, a class of entropy functionals were established to yield uniform $L^p$ estimates.

We conjecture that the bounded weak solution of semiconductor model is unique. Once $\mu_1\neq\mu_2$, then for weak solutions $u=(u_1,u_2),v=(v_1,v_2)$ of (\ref{equation})-(\ref{diffusiontwo}) and a test function $\varphi$, we have
\begin{equation}\label{deflectrelation}
\begin{aligned}
&\langle\partial_t(\mu_2u_1+\mu_1u_2),\varphi\rangle=\mu_1\mu_2\langle\Delta(u_1+u_2),\varphi\rangle,\\&\langle\partial_t(\mu_2v_1+\mu_1v_2),\varphi\rangle=\mu_1\mu_2\langle\Delta(v_1+v_2),\varphi\rangle.
\end{aligned}
\end{equation}

If $\mu_1=\mu_2$, by (\ref{deflectrelation}), we can show that $u_1+u_2=v_1+v_2$, a.e. $(x,t)\in\Omega\times(0,\infty)$. Then we apply the distance functional to show that $u=v$, a.e. $(x,t)\in\Omega\times(0,\infty)$.

If $\mu_1\neq\mu_2$, we cannot show that $\mu_2u_1+\mu_1u_2=\mu_2v_1+\mu_1v_2$ or $u_1+u_2=v_1+v_2$, a.e. $(x,t)\in\Omega\times(0,\infty)$ merely by the relation (\ref{deflectrelation}). Then we have difficulties in applying the distance functional to show that $u=v$, a.e. $(x,t)\in\Omega\times(0,\infty)$. The uniqueness discussion for semiconductor model is very difficult, and we are still working on it. 

For Shigesada-Kawasaki-Teramoto type population system, we also conjecture that once the initial value is bounded, then a bounded weak solution exists. 

The structure of diffusion matrix (\ref{sktdiffusion}) is rather complicated, that so far this conjecture is still being investigated. We have made lots of attempts, including the transformation technique, to show the population system has a bounded weak solution.
\\
\\

\textbf{Author contributions.} The author participated in the writing of the manuscript.
\\

\textbf{Funding.} The author acknowledges support from the National Natural Science
Foundation of China (NSFC), grant 12471206.
\\

\textbf{Data Availability.} No datasets were generated or analysed during the current study.
\\

\textbf{Competing interests} The authors declare no competing interests.
\\

%\subsection*{Acknowledgment}

\end{document}